\documentclass[12pt, a4paper]{article}
\usepackage[latin1]{inputenc}
\usepackage[english]{babel}
\usepackage{indentfirst}
\usepackage{amstext}
\usepackage{setspace}
\usepackage{mathtools}
\usepackage{amsfonts}
\usepackage{textcomp}
\usepackage{amssymb}
\usepackage{amscd}
\usepackage{epsf}
\usepackage{graphicx}
\usepackage{epsfig}

\newcommand {\demo}{\hskip -0.6cm{\bf Proof:  }}
\newcommand {\fim}{\hfill{$\square$}\vskip 1pc}

\newcommand{\s}{\sigma}

\newcommand {\N}{\mathbb{N}}

\newtheorem{theorem}{Theorem}[section]

\newtheorem{proposition}[theorem]{Proposition}

\begin{document}

\title{$(M + 1)$-step shift
spaces that are not conjugate to $M$-step shift spaces}
\maketitle
\begin{center}

{\large Daniel Gonçalves\footnote{This author is partially supported by CNPq.} and Danilo Royer}\\
\end{center}  
\vspace{8mm}

\abstract Recently Ott, Tomforde and Willis proposed a new approach for one sided shift spaces over infinite alphabets. In this new approach the conjugacy classes of shifts of finite type, edge shifts, and $M$-step shifts are distinct and the authors conjecture that for each $M\in \N\cup \{0\}$ there exist an $(M+1)$-step shift space that is not conjugate to any $M$-step shift. In this short paper we build a class of $(M+1)$-step shifts that are not conjugate to any $M$-step shift and hence show that their conjecture is correct.
\doublespace\newline\newline
MSC 2010: 37B10, 54H20

\section{Introduction}

In the past decades there have been numerous attempts to define and develop a theory of shift spaces over infinite alphabets (see \cite{OTW} for a comprehensive list), one of the main difficulties being the fact that when dealing with infinite alphabets the associated shift spaces are usually not compact. Recently, in \cite{OTW}, Ott, Tomforde and Willis proposed a new approach for one sided shift spaces over infinite alphabets, in which all shift spaces are compact. Among other key features of this new approach are applications to Cuntz-Krieger and Leavitt path algebras and graph C*-algebras, as well as the fact that many classical results of the theory of shifts over finite alphabets can be generalized. 

Although the approach in \cite{OTW} allows many results from classical symbolic dynamics to be generalized to the infinite alphabet case, there are also some important differences. In particular, in the finite alphabet case, it is well known that a shift space is a shift of finite type (i.e.,
described by a finite set of forbidden blocks) if and only if it is conjugate
to the edge shift coming from a finite graph with no sinks if and only if it
is an $M$-step shift (i.e., the shift is described by a set of forbidden blocks
all of length $(M + 1)$), see for example \cite{LM}. This is not the case in the infinite alphabet, where
the conjugacy classes of shifts of finite
type, edge shifts, and $M$-step shifts are distinct and hence have to be studied separately. 

In \cite{OTW} the authors show that $$\{0-\text{step shifts}\} \subseteq \{1-\text{step shifts}\}\subseteq \{2-\text{step shifts}\} \subseteq \ldots $$
and, at the end of Section 5, ask: 

\vspace{0.5pc}
For each $M\in \N\cup \{0\}$ does there exist an $(M + 1)$-step shift space that is not conjugate to any $M$-step shift?
\vspace{0.5pc}

Ott, Tomforde and Willis conjecture that the answer to the above question is yes and, in Section 3, we show that their conjecture is correct. In fact we will build a infinite family of $(M+1)$-step shift spaces that are not conjugate to any $M$-step shift space. In order to make the paper self contained, in Section 2 below we recall the relevant concepts regarding one sided infinite shifts over infinite alphabets as defined in \cite{OTW}.


\section{One sided shift spaces over infinite alphabets}

In this section we present a brief introduction to one sided shift spaces over infinite alphabets, as defined in \cite{OTW}. 

Given a set (alphabet) $A$, define a new symbol, say $\O$, which will indicate the {\em empty sequence} and let
$\Sigma_A^{fin}:=\{\O\}\cup\bigcup\limits_{k\in\N}A^k$ and $\Sigma_A^{inf}:=\{(x_i)_{i\in\N}:\ x_i\in A\ \forall i\in\N\}$.
We say that $\Sigma_A^{fin}$ is the set of all finite sequences over $A$, while $\Sigma_A^{inf}$ is the set of all infinite sequences over $A$.

The {\em full shift} over $A$ is the set $$\Sigma_A:=\left\{\begin{array}{lcl}\Sigma_A^{inf} & \text{ if} & |A|<\infty\\ \Sigma_A^{inf}\cup\Sigma_A^{fin} & \text{ if} & |A|=\infty.\end{array}\right.$$
Given $x\in\Sigma_A$ we define $l(x)$ as the length of $x$, which is equal to $k$ if $x\in A^k$ and equal to $\infty$ otherwise.

Next we recall the basis for the topology in $\Sigma_A$. Given $x=(x_1,\ldots,x_k)\in\Sigma_A^{fin}$, $x\neq\O$, and a finite set $F\subset A$, we define a {\em generalized cylinder set} of $\Sigma_A$ as
$Z(x,F):=\{y\in\Sigma_A:\ y_i=x_i\ \forall i=1\ldots,k,\ y_{k+1}\notin F\}$. If $ x = \O $ then we let $Z(\O,F)= \{y\in\Sigma_A:\ y_1\notin F\}$.
														
If $F=\emptyset$ we will denote $Z(x,F)$ just by $Z(x)$. Notice that $Z(x)$ coincides with the usual cylinders which form the basis of the product topology in $\Sigma_A^{inf}$, while $Z(\O,F)$ is a neighborhood of the empty sequence $\O$. It is proven in \cite{OTW} that the collection of generalized cylinder sets defined above forms a basis for a topology in $\Sigma_A$ consisting of compact open sets. The full shift space over $A$ is defined as $\Sigma_A$ with the topology above and in \cite{OTW} it is proved that $\Sigma_A$ is zero dimensional (that is, it has a basis of clopen sets) and compact.

The {\em shift map} $\s:\Sigma_A\to\Sigma_A$ is defined as follows:
$$\s(x)=\left\{\begin{array}{lcl} \O & \text{ if} & x=\O \text{ or } x\in A^1\\
                      (x_{i+1})_{i=1}^{k-1} & \text{ if} & x=(x_i)_{i=1}^k\in A^k\\
                 (x_{i+1})_{i\in\N} & \text{ if} & x=(x_i)_{i\in\N}\in \Sigma_A^{inf}.\end{array}\right.$$ 
We remark that, if $A$ is infinite, than $\s$ is not continuous at the empty sequence $\O$ (Proposition 2.3 in \cite{OTW}).


We will be particularly interested in special subsets of $\Sigma_A$, which are named shift spaces.
We say that $\Lambda\subseteq\Sigma_A$ is a {\em shift space} over $A$ if the following three properties hold:
\begin{enumerate}
\item[1] $\Lambda$ is closed with respect to the topology of $\Sigma_A$;
\item[2] $\Lambda$ is invariant under the shift map, that is, $\s(\Lambda)\subseteq\Lambda$;
\item[3] $\Lambda$ satisfies the `infinite-extension property', that is, for all $x \in \Lambda$, with $l(x)< \infty$, the set $\{a \in A: xay \in \Lambda \text{ for some } y \in \Sigma_A\}$ is infinite.
\end{enumerate}

Notice that properties 1 and 2 above are exactly the same ones that define a shift space over a finite alphabet. Property 3 ensures that there always exists infinite sequences in a non-empty shift space. It is proven in Proposition~3.8 of \cite{OTW} that $\Lambda^{inf}$ is dense in $\Lambda$.

In symbolic dynamics, an equivalent way to define a shift space is via a set of forbidden words: Let $\mathbf{F}\subset \bigcup\limits_{k\geq 1}A^k$, 

$\begin{array}{ll} X_{\mathbf{F}}^{inf} = & \{ x \in \Sigma_A^{inf}: \text{ no subblock of $x$ is in $\mathbf{F}$} \}, \text{ and } \\ 
X_{\mathbf{F}}^{fin} = & \{ x \in \Sigma_A^{fin}: \text{ there are infinitely many } a \in A \\
& \text{ for which there exists } y\in \Sigma_A^{inf} \text{  such that } xay\in  X_{\mathbf{F}}^{inf}  \},
\end{array}$
where a subblock of $x$ means an element $u\in \Sigma_A^{fin}$ such that $x=yuz$ for some $y\in \Sigma_A^{fin}$ and some $z \in \Sigma_A$. Now, define $ X_{\mathbf{F}}$ as $X_{\mathbf{F}}^{inf} \cup X_{\mathbf{F}}^{fin}$. It is proven in \cite{OTW}, theorem 3.16, that $\Lambda$ is a shift space iff $\Lambda = X_{\mathbf{F}}$ for some $\mathbf{F}\subset \bigcup\limits_{k\geq 1}A^k$.

Following \cite{OTW}, we say that $\Lambda$ is a shift of finite type, SFT, if we can take $\mathbf{F}$ having only finitely elements and an $M$-step shift if $l(x)=M+1$ for all $x\in \mathbf{F}$. 


Finally recall that a map $\phi: \Lambda \rightarrow Y$, where $\Lambda\subseteq \Sigma_A$, $Y\subseteq \Sigma_B$ are shift spaces over alphabets $A$ and $B$ respectively, is a conjugacy if, and only if, it is bijective, continuous, shift commuting and $l(\phi(x)) = l(x)$ for all $x \in \Lambda$.

We can now proceed to the main section or our work.

\section{A family of $(M + 1)$-step shifts that are not conjugate to any $M$-step shift}

In order to construct the relevant $M$-step shifts we need to view the forbidden words in a shift $X_{\mathbf{F}}$ as the inverse image of a function. More precisely, given $\mathbf{F}\subseteq \bigcup\limits_{k\geq 1}A^k$ we can see $\mathbf{F}$ as the inverse image $f^{-1}(0)$, where $f:\bigcup\limits_{k\geq 1}A^k\rightarrow \{0,1\}$ is given by $f(x)=0$, if $x\in \mathbf{F}$ and $f(x)=1$, if $x\notin \mathbf{F}$. 

Given $f:\bigcup\limits_{k\geq 1}A^k\rightarrow \{0,1\}$, let 

$\begin{array}{ll} X_f^{inf}= &\{x\in \Sigma_A^{inf} \text{ such that } f(x_1,...,x_n)=1 \\ & \text{ for each subblock } x_1...x_n \text{ of } x\}, \\

X_f^{fin}=& \{x\in \Sigma_A^{fin} \text{ such that there are infinitely many }a\in A \\ &   \text{ for which there exists } y\in \Sigma_A^{inf} \text{ such that } xay \in X_f^{inf}\}
\end{array}$
and define
$X_f=X_f^{inf}\cup X_f^{fin}$.

With the above definitions one can readily see that any shift can be written as $X_f$ for a suitable function $f$. 
In particular, a shift space $\Lambda\subseteq \Sigma_A$ is an $M$-step shift if, and only if, $\Lambda=X_f$ for some $f:\bigcup\limits_{k\geq 1}A^k\rightarrow \{0,1\}$ induced by an auxiliary function $\tilde{f}:A^{M+1}\rightarrow \{0,1\}$ as follows:
$$
f(x)=\left\{\begin{array}{ll}
 1  & \text {if} |x|\leq M\\
\prod\limits_{i=1}^{|x|-M}\tilde{f}(x_i,...,x_{M+i}) & \text{if} |x|\geq M+1 .
\end{array}\right .$$

A map $f$ defined from an auxiliary function $\tilde{f}:A^{M+1}\rightarrow \{0,1\}$ as above will be called an {\it $M$-step shift map induced by $\tilde{f}$.}

We can now proceed with the construction of the shift spaces. Let $A$ be an infinite alphabet and $M\in \N$ with $M\geq 1$. In light of the discussion above, all we need to do to construct an $(M+1)$-step shift is to define an auxiliary function $\tilde{f}:A^{M+2} \rightarrow \{0,1\}$. In particular the relevant auxiliary maps $\tilde{f}:A^{M+2} \rightarrow \{0,1\}$ satisfy the following two conditions:

\begin{itemize}
\item[(1)] there exists a subset $\{x_1,...,x_M\}\subseteq A$ such that there are infinitely many $x_{(M+1)}\in A$ for which 

$\begin{array}{l} 1=\tilde{f}(x_1,x_2,...x_M,x_{(M+1)},x_1)=\tilde{f}(x_2,x_3,...,x_{(M+1)},x_1,x_2)=\cdots \\
\cdots=\tilde{f}(x_M,x_{(M+1)},x_1,...,x_M)=\tilde{f}(x_{(M+1)},x_1,x_2,...,x_{(M+1)}). 
\end{array}$

\item[(2)] for each $\{x_1,...,x_{(M+1)}\}\subseteq A$ there are only finitely $x_{(M+2)}\in A$ such that $\tilde{f}(x_1,...,x_{(M+1)},x_{(M+2)})=1$.
\end{itemize}


\begin{proposition}\label{prop1} Suppose that $M\geq 1$, let $\tilde{f}:A^{M+2}\rightarrow \{0,1\}$ be a map satisfying conditions $(1)$ and $(2)$ above, let $f$ be the $(M+1)$-step shift map induced by $\tilde{f}$ and $X_f$ be the associated $(M+1)$-step shift. Then $X_f$ is not conjugate to any $M$-step shift space $Y$.
\end{proposition}

\demo We will give a proof by contradiction. Suppose that there exists a conjugacy $\Phi:X_f\rightarrow Y_g$, where $Y_g$ is an $M$-step shift over the alphabet $B$, associated with an $M$-step shift map $g$ induced by an auxiliary map $\tilde{g}:B^{(M+1)}\rightarrow \{0,1\}$. 

Let $\{x_1,x_2,...,x_M\}\subseteq A$ be as in condition (1) and let $\{x_{(M+1)}^j:j\in \N\}$ be an infinite set of distinct elements such that the equality in condition (1) holds. We will use these elements to construct a sequence of periodic elements in $X^{inf}_f$ converging to an element in $X^{fin}_f$. Precisely, for each $j \in \N$, let 
$$\xi^j=x_1x_2x_3...x_Mx_{(M+1)}^jx_1x_2x_3...x_Mx_{(M+1)}^j\cdots.$$ 
Notice that $\xi^j \in X^{inf}_f$ and the sequence $(\xi^j)_{j\in \N}$ converges to $\xi=x_1...x_M\in X_f$.

We now look at the image of the sequence above in $Y_g$. For clarity, let $\Phi(\xi^j)=y^j$ and $\Phi(\xi)=y$. Since $\Phi$ is a conjugacy, $(y^j)_{j\in N}$ converges to $y$, $l(y^j)=\infty$, $l(y)=M$ and hence we can write $y=y_1\cdots y_M$ and $y^j=y_1^jy_2^jy_3^j\cdots$, where $y_i, y_i^j\in B$ for each $i,j$. Furthermore, since $\xi^j$ has period $(M+1)$, $y^j$ also has period $(M+1)$ (that is $\sigma^{M+1}(y^j)=y^j$) for all $j\in \N$.

Now, since $(y^j)_{j\in \N}$ converges to $y$, there exist a subsequence $(y^{j_k})_{k\in \N}$ such that the beginning of each $y^{j_k}$ is equal to $y$, that is $$y^{j_k}=y(y_{(M+1)}^{j_k})(y_{(M+2)}^{j_k})(y_{(M+3)}^{j_k})\cdots$$ and the elements $y_{(M+1)}^{j_k}$ are all pairwise distinct. Moreover, since each $y^{j_k}$ has period $(M+1)$, we can write, for each $j^k$,
$$y^{j_k}=y(y_{(M+1)}^{j_k})y(y_{(M+1)}^{j_k})y(y_{(M+1)}^{j_k})\cdots,$$ and since $y^{j_k}\in Y_g^{inf}$ we have that, for each $j_k$, $$1=\tilde{g}(y_1...y_My_{(M+1)}^{j_k})=\tilde{g}(y_2...y_My_{(M+1)}^{j_k}y_1)=\tilde{g}(y_3...y_My_{(M+1)}^{j_k}y_1y_2)\cdots$$
$$\cdots=\tilde{g}(y_M y_{(M+1)}^{j_k}y_1...y_{(M-1)})=\tilde{g}(y_{(M+1)}^{j_k}y_1...y_M).$$

Next, let $j_{k_0}\in \N$ be fixed and, for all $j_k$, define $$z^{j_k}=y_1...y_My_{(M+1)}^{j_{k_0}}y^{j_k}.$$
Using the above equalities involving $\tilde{g}$, and the definition of the $M$-step shift map $g$ induced by $\tilde{g}$, one can easily check that $g(r)=1$ for each subblock $r$ of $z^{j_k}$ and hence $z^{j_k}\in Y_g^{inf}$ for all $j_k$. 

Finally, to reach a contradiction, notice that $(z^{j_k})_{k\in \N}$ converges to $y_1...y_My_{(M+1)}^{j_{k_0}}y_1...y_M$, which has length $(2M+1)$. So, $Y_g$ contains an element of length $(2M+1)$ and, since $\Phi$ is a conjugacy, $X_f$ must also contain an element of length $(2M+1)$. However, it follows from condition (2) satisfied by $\tilde{f}$ that $X_f$ contains no elements of length $(M+2)$ and hence it cannot contain elements of length $(2M+1)$.

\fim


With the above proposition we can completely answer Ott, Tomforde and Willis's question. We do this below. 

\begin{theorem}
For every $M\in \N\cup \{0\}$ there exists a one-sided shift
space over an infinite alphabet that is an $(M + 1)$-step shift space and is
not conjugate to an M-step shift space".
\end{theorem}
\demo
First suppose $M\geq 1$.  Let $\tilde{f}:A^{M+2}\rightarrow \{0,1\}$ be defined by $$\tilde{f}(x_1,...,x_{(M+2)})= \begin{cases} 1  & \text{ if }x_1=x_{(M+2)} \\0& \text{ if } x_1\neq x_{(M+2)}.
\end{cases}$$ Then $\tilde{f}$ is a function satisfying the hypothesis of Proposition~\ref{prop1} and the theorem follows by direct aplication of Proposition~\ref{prop1}

Next suppose $M=0$. Let $X_f$ be the shift space associated with the $1$-step shift map $f$ induced by an auxiliary map $\tilde{f}:A^2\rightarrow \{0,1\}$ satisfying the following two conditions: 
\begin{enumerate}
\item[(i)]for each $x\in A$ there is an $y\in A$ such that $\tilde{f}(x,y)=1$;
\item[ (ii)] the set $\{x\in A: \tilde{f}(x,y)=1 \text{ for a infinite number of } y\in A\}$ is finite and non-empty. 
\end{enumerate}

An example of such a function is, for $x_0\in A$ fixed, $\tilde{f}(x,y)=1$ if $x\neq x_0$ and $y=x$, $\tilde{f}(x,y)=1$ if $x=x_0$ and $\tilde{f}(x,y)=0$ otherwise.

Notice that $X_f$ contains only finitely many elements with length 1, namely the elements $x\in A$ for which there are infinitely many $y\in A$ such that $\tilde{f}(x,y)=1$. 

If we now suppose that $X_f$ is conjugate to a 0-step shift space $Y$ (which is a full shift over some alphabet $B$) then we have that $Y$ contains some element of length $1$, and so $B$ is infinite. But this implies that $Y$ contains infinitely many elements of length 1 and hence so does $X_f$, which is a contradiction. We conclude that $X_f$ is not conjugate to any 0-step shift space. 
\fim


\end{document}